\documentclass[12pt]{amsart}
\usepackage{amssymb}

\newcommand{\prop}[1]{\noindent\textbf{#1}.}
\newtheorem*{theorem}{Theorem}

\newcommand{\R}{\mathbb R}
\newcommand{\A}{\mathfrak A}
\newcommand{\supp}{\operatorname{supp}}
\newcommand{\osc}{\operatorname{osc}}
\newcommand{\ci}[1]{_{{}_{\scriptstyle{#1}}}}
\newcommand{\wt}{\widetilde}
\newcommand{\e}{\varepsilon}
\newcommand{\f}{\varphi}
\renewcommand{\le}{\leqslant}
\renewcommand{\ge}{\geqslant}

\begin{document}

\title[The $s$-Riesz transform of an $s$-dimensional measure]
{The $s$-Riesz transform of an $s$-dimensional measure in $\R^2$ is unbounded for $1<s<2$}

\author{Vladimir Eiderman}
\address{Vladimir Eiderman,  Department of  Mathematics, University of Wisconsin-Madison, Madison, WI}
\email{eiderman@math.wisc.edu}
\author{Fedor Nazarov}
\address{Fedor Nazarov, Department of Mathematics, Kent State University, Kent, OH}
\email{nazarov@math.kent.edu}
\author{Alexander Volberg}
\address{Alexander Volberg, Department of  Mathematics, Michigan State University, East Lansing, MI}
\email{volberg@math.msu.edu}
\thanks{Work of F.~Nazarov and  A.~Volberg is supported  by the National Science Foundation under the grant  DMS-0758552. 
}
\begin{abstract}
In this paper, we prove that for $s\in(1,2)$ there exists no totally lower irregular finite positive Borel measure $\mu$ in $\R^2$ with\break $\mathcal H^s(\supp\mu)<+\infty$ such that 
$\|R\mu\|\ci{L^\infty(m_2)}<+\infty$, where $R\mu=\mu\ast\frac{x}{|x|^{s+1}}$ and $m_2$ is the Lebesgue measure in $\R^2$. Combined with known results of Prat and  Vihtil\"a, this shows that for  any non-integer $s\in(0,2)$ and any finite positive Borel measure in $\R^2$ with $\mathcal H^s(\supp\mu)<+\infty$, we have  $\|R\mu\|\ci{L^\infty(m_2)}=\infty$.
\end{abstract}
\maketitle

\section{Introduction}\label{introd}

Let $\mu$ be a finite strictly positive Borel measure on the plane $\R^2$.
We will say that $\mu$ is $s$-dimensional if $\mathcal H^s(\supp\mu)<+\infty$ where
$\mathcal H^s$ is the $s$-dimensional Hausdorff measure. Another way to state it is that 
there exists some positive $H<+\infty$ such that for every $r>0$, one can find a (countable)
sequence of disks $D_i=D(c_i,r_i)$ with centers $c_i$ and radii $r_i$ such that
 $r_i<r$ for all $i$, $\sum_i r_i^s\le H$, and 
$\mu(\R^2\setminus\cup_i D_i)=0$. 

An $s$-dimensional measure $\mu$ is called totally lower irregular if 
$$
\liminf_{r\to 0+}r^{-s}\mu(D(x,r))=0\qquad\text{for $\mu$-a.e. }x\in\R^2\,.
$$  
If $\nu$ is a finite (signed) measure on $\R^2$, its ($s$-dimensional) Riesz transform
$R\nu$ is defined by 
$$
(R\nu)(x)=\int_{\R^2}\frac{x-y}{|x-y|^{s+1}}\,d\nu(y)\,.
$$
If $0<s<2$, the integral in this definition converges absolutely almost everywhere with respect
to the 2-dimensional Lebesgue measure $m_2$ on $\R^2$. 
If, in addition to being finite, $\nu$ has bounded density with respect to $m_2$, the integral
converges everywhere and is a continuous function on the plane that tends to $0$ at infinity.

We will say that $R\nu$ is bounded if $\|R\nu\|\ci{L^\infty(m_2)}<+\infty$.

Our goal is to complete the proof of the following theorem.

\begin{theorem} Let $s\in(0,2)\setminus\{1\}$ and let $\mu$ be a strictly positive finite Borel measure in $\R^2$ such that $\mathcal H^s(\supp\mu)<+\infty$. Then
$\|R\mu\|\ci{L^\infty(m_2)}=\infty$.
\end{theorem}

It is easy to see that for $s=1,2$ this statement is incorrect. Indeed, for any nonnegative $\f\in C_0^\infty(\R^2)$, the measures $\mu=\f\mathcal H^1|_L$, where $L$ is some line in $\R^2$, and $\mu=\f m_2$ give counterexamples for $s=1$ and $s=2$ respectively. 

For non-integer $s\in(0,2)$, the theorem has been known in the following cases.

For $0<s<1$, it has been proved by Prat \cite{P} using Melnikov's curvature techniques introduced  in \cite{M}. Unfortunately this tool is ``cruelly missing'' (by the expression of Guy David) for $s>1$, because the natural analog of the squared Menger curvature can be negative.

For Riesz transforms in $\R^d$ corresponding to non-integer $s\in(0,d)$, the unboundedness of $R\mu$ was established by Vihtil\"a \cite{Vi} under the additional assumption that the lower $s$-density of $\mu$ is positive for $\mu$-almost all $x\in\supp\mu$. The main tool in \cite{Vi} is the concept of the tangent measure. This method also fails  in the general case because without any assumptions on lower density, the tangent measure may lose the property of being $s$-dimensional. On the other hand, \cite{Vi} gives more than is formally claimed there. The same argument (but if one adds some non-homogeneous Harmonic Analysis consideration like in \cite{VoRMI} for example)  yields the desired assertion for any finite measure $\mu$ such that $\mu\{x:\liminf_{r\to 0+}r^{-s}\mu(D(x,r))>0\}>0$. Thus, to finish the proof of the theorem, it is enough to consider the case of $s$-dimensional totally irregular measures, which is exactly what we will do in the current paper. This requires introducing several new techniques, which, we hope, may be of independent interest.

Note that our theorem, as well as the results of Prat and Vihtil\"a, apply to arbitrary $s$-dimensional measures. When $\mu$ is supported on a Cantor set of certain type in $\R^d$, the unboundedness of its Riesz transform follows immediately from explicit bounds for Calder\'on-Zygmund capacities of Cantor sets in \cite{MT}, \cite{T1},\cite{EV}, and other similar papers.

It is also worth mentioning that de Villa and Tolsa \cite{RVXT} proved that the Riesz transform of  an $s$-dimensional measure in $\R^d$ cannot have principal values for non-integer $s$.

\section{Definitions and notation}\label{defnot}

The operator $R$ returns a vector-valued function and is often written as $(R_1,R_2)$ where 
$R_j\nu$ is the $j$-th coordinate of $R\nu$ ($j=1,2$). We shall denote by $R^*$ the formal 
adjoint of $R$ that acts on vector-valued finite measures $\eta$ by the rule 
$R^*\eta=-\sum_j R_j\eta_j$ where $\eta_j$ are the ``coordinate measures'' of $\eta$. The
identity
$$
\int_{\R^2}\langle R\nu,d\eta\rangle=\int_{\R^2}R^*\eta\,d\nu
$$
holds every time when at least one of the finite measures involved has bounded density 
with respect to $m_2$ (here $\langle\cdot,\cdot\rangle$ denotes the scalar product in $\R^2$).

By $C$ with or without an index we shall denote a (large)
positive constant that may depend only on $s$. This constant
may change from line to line if it has no index. The indexed constants are fixed throughout
the paper and the convention is that $C_j$ can be chosen as soon as all $C_i$ with $i<j$
are known.

For reader's convenience, we will list a few symbols that will occur rather frequently. 
\medskip

$s$ a number in $(1,2)$;

$D(x,r)$ the disk of radius $r$ centered at $x$;

$\A=\{2,4,8,16,\dots\}$ the set of positive integer powers of $2$;

$\nu,\eta$ generic measures (possibly signed or even vector-valued);

$N,\e,M,\delta$ positive parameters to be chosen in this order. $N$ and $M$ are large,
$\e$ and $\delta$ are small;

$\mu$ the totally lower irregular $s$-dimensional measure with bounded Riesz transform
(i.e., the measure whose non-existence we want to prove);

$m$ one half of the total mass of $\mu$;

$H$ twice the Hausdorff measure of the support of $\mu$;

$\mu'$ the part of $\mu$ obtained by dropping everything supported outside
the lowest level of the Cantor construction;

$\wt\mu$ the mollified $\mu'$ with smooth density consisting of small caps supported on 
$\wt\Omega_j$;

$\wt\Omega_j$ the disk of radius $\e\rho_j$ contained in $\Omega_j$;

$\Omega_j$ the $\e\rho_j$-neighborhood of $\wt B_j$;

$\wt B_j=(1-3\e)B_j\setminus \cup_{i<j}B_i$;

$B_j$ the disks in the bottom cover;

$T_j$ the disks in the top cover;

$\wt T_j=T_j\setminus \cup_{i<j}T_i$;

$\psi$ the vector-valued function associated with the top cover;

$\Psi$ the majorant of $|\psi|$;

$R$ the Riesz transform;

$R^*$ the adjoint Riesz transform;

$R^\sharp$ the maximal Riesz transform;

$(R^*)^\sharp$ the maximal adjoint Riesz transform;

$U$ the $(s-1)$-dimensional Newton potential;

$V$ a smooth convex version of $\min(|x|^2,|x|)$.
\medskip

We shall also assume that the reader is familiar with the basic theory of singular integral
operators in non-homogeneous spaces.     

\section{Elementary properties of the Riesz transform}\label{elepro}

We shall use the following standard facts without any special references.

\prop{Translation invariance and scaling} If $f\in L^1(m_2)$, then
$$
[R(f(\tfrac{\cdot-c}{r})m_2)](x)=r^{2-s}[R(fm_2)](\tfrac{x-c}{r})\,,\quad x\in\R^2.
$$
\prop{Action on the Fourier side}
$$
\widehat{R_j(fm_2)}(\xi)=i\sigma\frac{\xi_j}{|\xi|^{3-s}}\widehat f(\xi)\,,
$$
where $\sigma\ne 0$ is some real constant. 

More precisely, if $f,g,\widehat f,\widehat g\in L^{1}(m_2)\cap L^\infty(m_2)$, then 
$$
\int_{\R^2}[R_j(fm_2)]\bar g\,dm_2=i\sigma\int_{\R^2}\frac{\xi_j}{|\xi|^{3-s}}\widehat f(\xi)
\overline{\widehat g(\xi)}\,dm_2(\xi)\,.
$$

\prop{The $L^\infty$ bound} If $\supp f$ is contained in a disk of radius $r$, then 
$$
\|R(fm_2)\|\ci{L^\infty(m_2)}\le C_1\|f\|_{L^\infty(m_2)}r^{2-s}\,.
$$

\prop{Relation to the Newton potential} Let 
$U\nu(x)=-\frac{1}{s-1}\int_{\R^2}\frac{d\nu(y)}{|x-y|^{s-1}}$. Then 
$$
R_j\nu=\frac{\partial}{\partial x_j}U\nu\,.
$$
If $\nu=fm_2$ with smooth compactly supported $f$, we can pass the derivative to $f$ and write
$$
R_j\nu=U(\tfrac{\partial f}{\partial x_j}m_2)\,.
$$

\section{The representation of the standard cap}

Let $\f_\circ$ be any positive Schwartz function that is at least $1$ on the unit disk
centered at the origin. Define the vector field $\psi_\circ$ by
$$
\widehat\psi_\circ(\xi)=i\sigma^{-1}\xi|\xi|^{1-s}\widehat\f_\circ(\xi)\,.
$$
We claim that 
$$
|\psi_\circ(x)| \le\frac{C_2}{(1+|x|)^{4-s}}\qquad\text{and}\qquad R^*(\psi_\circ m_2)=\f_\circ\,.
$$
The second claim follows from the first at once if we check the action of both sides on nice 
test-functions and pass to the Fourier side (which is justified because 
$\psi_\circ,\widehat\psi_\circ\in L^1(m_2)\cap L^\infty(m_2)$). The first claim is a standard 
exercise in elementary Fourier analysis left to the reader.

\section{The growth bound and its implications}\label{grobou}

Let $\mu$ be a finite positive measure satisfying $\|R\mu\|\ci{L^\infty(m_2)}\le 1$. Take a disk
$D=D(c,r)$ and write
\begin{multline*}
\mu(D)\le\int_{\R^2}\f_\circ(\tfrac{\cdot-c}{r})\,d\mu=
\int_{\R^2}R^*[r^{s-2}\psi_\circ(\tfrac{\cdot-c}{r})m_2]\,d\mu
\\
=r^s\int_{\R^2}\langle R\mu, r^{-2}\psi_\circ(\tfrac{\cdot-c}{r})\rangle\,dm_2
\le r^s\|R\mu\|\ci{L^\infty(m_2)}\|\psi_\circ\|\ci{L^1(m_2)}\le C_3 r^s\,.
\end{multline*}
This a priori growth bound combined with the assumption 
$$
\|R\mu\|\ci{L^\infty(m_2)}\le 1
$$
allows one to apply to the measure $\mu$ the whole non-homogeneous singular integral machinery (see, e.g.. \cite{NTV-IMRN98}, \cite{NTV-Acta})
and to conclude that the maximal singular operators
$f\mapsto R^\sharp(f\mu)$ and $g\mapsto (R^*)^\sharp(g\mu)$ are bounded in $L^2(\mu)$ 
with norms not exceeding $C_4$. Here $R^\sharp$ is the maximal Riesz transform defined by 
$$
(R^\sharp\nu)(x)=\sup_{D:x\in D}\Bigl|\int_{\R^2\setminus 2D}\frac{x-y}{|x-y|^{s+1}}\,d\nu(y)
\Bigr|\,,
$$
where the supremum is taken over all disks $D$ containing $x$, and $2D$ stands for the
disk with the same center as $D$ but of twice larger radius. The operator $(R^*)^\sharp$ is defined in a similar way.

Note that, unlike the initial assumption $\|R\mu\|\ci{L^\infty(m_2)}\le 1$, the growth bound
and the operator norm condition are preserved if we drop any part of the measure $\mu$.
In what follows, we will rely on these two conditions only and never use the $L^\infty$ bound itself.

From now on, $\mu$ will be a fixed finite measure of total mass $2m$, supported on a set of ($s$-dimensional) Hausdorff measure $H/2$, and satisfying the 
growth bound and the operator norm condition above. Note that the growth bound implies that
we automatically have $m\le C_3 H$.

\section{The good old Cantor set argument}\label{gooold}
 
The main motivation for our construction is the following well-known argument 
for Frostman measures on sparse Cantor squares. Assume that we have a sparse Cantor
square $K$ of dimension $s$ 
on the plane in which the squares of each generation are separated by
distances much larger than their diameters.

For $x\in K$, 
let $K^{(n)}(x)$ be the square of the $n$-th generation 
containing $x$. Let $\nu=\mathcal H^s|\ci K$. Define
$$
R^{(n)}\nu(x)=\int_{K^{(n)}(x)\setminus K^{(n+1)}(x)}\frac{x-y}{|x-y|^{s+1}}\,d\nu(y)\,.
$$
Then $R^\sharp\nu$ dominates every partial sum $\sum_{n=0}^{N-1}R^{(n)}\nu$.

The key observation is that the norms $\|R^{(n)}\nu\|\ci{L^2(\nu)}$ are uniformly bounded from below in $L^2(\nu)$ because for every $x$, the differences $x-y$ point pretty much in the same
direction when $y\in K^{(n)}(x)\setminus K^{(n+1)}(x)$ and the kernel blow-up
near the diagonal is perfectly balanced with the decay of the measure. On the
other hand, the oscillation $\osc_{K_j^{(n+1)}}R^{(n)}\nu$ is very small for
every Cantor square $K_j^{(n+1)}$ of the $(n+1)$-st generation and we also have the
cancellation property 
\begin{multline*}
\int_{K_j^{(n+1)}}[R^{(n+1)}\nu]\,d\nu= \\
\iint_{x,y\in K_j^{(n+1)}, K^{(n+2)}(x)\ne K^{(n+2)}(y)}\frac{x-y}{|x-y|^{s+1}}
\,d\nu(x)\,d\nu(y)   =0\,.
\end{multline*}
Together they
imply that the functions $R^{(n)}\nu$ are almost orthogonal in $L^2(\nu)$, so
$$
\int_{\R^2}\Bigl|\sum_{n=0}^{N-1}R^{(n)}\nu\Bigr|^2\,d\nu
\approx \sum_{n=0}^{N-1}\int_{\R^2}|R^{(n)}\nu|^2\,d\nu\approx N\,,
$$
and we can conclude that $R^\sharp$ is unbounded in $L^2(\nu)$.

We will use this simple argument as a guideline. The difficulty is that an arbitrary
$s$-dimensional set has no a priori Cantor type structure and an attempt to introduce
it using the standard dyadic scales encounters severe difficulties with both almost
orthogonality and the lower bounds for $R^{(n)}\mu$. We will use a slightly different
partition that gives the almost orthogonality for free in the case of totally lower
irregular measures. Still, we will have to fight hard for the lower bounds.

\section{The top cover and the associated $\Psi$-function}\label{topcov} 

Fix $N\in\mathbb N$, $\e>0$, $M>1$, $\delta>0$ to be chosen in this order.
The reader should think of $N,M$ as of very large parameters and of $\e,\delta$ as
of very small ones. Choose some $r^*>0$. We start with choosing a finite sequence of disks 
$T_j=D(c_j,r_j)$ such that $r_j\le r^*$, $\sum_j r_j^s\le H$, and $\mu(\R^2\setminus\cup T_j)<\e m$.
Without loss of generality, we may also assume that the union of the boundaries of
$T_j$ has zero $\mu$-measure. Put $\wt T_j=T_j\setminus \cup_{i<j}T_i$.

Define
$$
\psi=\sum_j\mu(\wt T_j)r_j^{-2}\psi_\circ(\tfrac{\cdot-c_j}{r_j})
$$
and 
$$
\Psi\ci A=\sum_j\frac{\mu(\wt T_j)}{\pi A^2r_j^2}\,\chi\ci{AT_j}\,,\qquad 
\Psi=\sum_{A\in\A}A^{s-2}\Psi\ci A\,,
$$
where $\A=\{2^k:k\in\mathbb N\}$ and $\chi\ci E$ is the characteristic function
of the set $E$.

Note that the pointwise bound for $\psi_\circ$ implies that $|\psi|\le C_5\Psi$. Also observe
that 
$$
R^*(\psi m_2)=\sum_j \frac{\mu(\wt T_j)}{r_j^s}\,\f_\circ(\tfrac{\cdot-c_j}{r_j})\ge
\sum_j \frac{\mu(\wt T_j)}{r_j^s}\,\chi\ci{\wt T_j}\,.
$$
Let $\nu$ be any finite positive measure supported on $\cup_jT_j$ and 
satisfying $\nu(\wt T_j)\le 2\mu(\wt T_j)$, $\nu(\R^2)\ge m$. Write
\begin{multline*}
C_5\int_{\R^2}|R\nu|\Psi\,dm_2\ge \int_{\R^2}\langle R\nu,\psi\rangle\,dm_2
=\int_{\R^2}R^*(\psi m_2)\,d\nu
\\
\ge\sum_j\frac{\mu(\wt T_j)}{r_j^s}\nu(\wt T_j)
\ge\frac 12\sum_j\frac{\nu(\wt T_j)^2}{r_j^s}
\\
\ge\frac 12\Bigl(\sum_j\nu(\wt T_j)\Bigr)^2\Bigl(\sum_j r_j^s\Bigr)^{-1}\ge \frac{m^2}{2H}\,.
\end{multline*}
On the other hand, we, clearly, have
$$
\int_{\R^2}\Psi\ci A\,dm_2=\sum_j\mu(\wt T_j)\le 2m\qquad\text{whence}\qquad 
\int_{\R^2}\Psi\,dm_2\le C_6m\,.
$$

\section{The function $V$}\label{functV}

Consider any $C^\infty$ function $v$ on $[0,+\infty)$ such that $v(0)=v'(0)=0$
and $v''$ is a non-increasing function that is identically $2$ on $[0,1]$ and
identically $0$ on $[2,+\infty)$. The function $v(t)$ is increasing, convex, 
equals $t^2$ on $[0,1]$, satisfies the inequalities $\min(t,t^2)\le v(t)\le t^2$ and $v'\le 4$
for all $t\ge 0$. Also we have $v'(t)=\int_0^tv''(\tau)\,d\tau\ge tv''(t)$, that is, $(tv')'\le(2v)'$. Hence, $tv'\le2v$. Integrating the inequality 
$\frac{v'(t)}{v(t)}\le\frac 2t$ from $t$ to $at$, $a\ge1$, we get $v(at)\le a^2v(t)$. Moreover, we have $v'(t)^2\le\frac{4v^2(t)}{t^2}\le 4 v(t)$.

Define $V(x)=v(|x|)$. In what follows, we will need a good lower bound for the integral
$\int_{\R^2}V(R\nu)\Psi\,dm_2$ under the same assumptions on $\nu$ as in the 
previous section. Put $I=\int_{\R^2}\Psi\,dm_2$ and apply
Jensen's inequality to get
\begin{multline*}
\int_{\R^2}V(R\nu)\Psi\,dm_2\ge Iv\left(I^{-1}\int_{\R^2}|R\nu|\Psi\,dm_2\right)\ge
Iv\left(I^{-1}\frac{m^2}{2C_5 H}\right)
\\
\ge\min\biggl(\frac{m^2}{2C_5 H},\frac{m^4}{4C_5^2H^2 I}\biggr)\ge
\min\biggl(\frac{m^2}{2C_5 H},\frac{m^3}{4C_5^2C_6 H^2 }\biggr)\ge C_7^{-1}\frac{m^3}{H^2}
\end{multline*}
(we used that $v(t)\ge\min(t,t^2)$, $\int_{\R^2} \Psi \,d m_2 \le C_6m$, and $m\le C_3H$ here).

\section{The Marcinkiewicz $g$-function}\label{marcin}

For $A\ge 2$, define 
$$
g\ci A=A^{-s}\sum_j\frac{\mu(\wt T_j)}{r_j^s}\,\chi\ci{AT_j}\,.
$$
We claim that $\int_{\R^2} g_A^2\,d\mu\le C_8m$. Indeed, let $f$ be any (positive) function
with $\|f\|\ci{L^2(\mu)}=1$. Then
$$
\int_{\R^2} g_A f\,d\mu=\sum_j \mu(\wt T_j)\frac{1}{(Ar_j)^s}\int_{AT_j}f\,d\mu
\le 3^sC_3\sum_j \mu(\wt T_j)\frac{1}{\mu(3AT_j)}\int_{AT_j}f\,d\mu
$$ 
because $\mu(3AT_j)\le C_3(3Ar_j)^s$. But the normalized integral factor is dominated by 
the non-homogeneous Hardy-Littlewood maximal function 
$$
\mathcal Mf(x)=\sup_{D:\,x\in D}\frac 1{\mu(3D)}\int_{D}f\,d\mu
$$ 
on $\wt T_j$. Thus, the last sum does not exceed 
$
\int_{\R^2} \mathcal Mf \,d\mu\le C\sqrt m\|\mathcal Mf\|\ci{L^2(\mu)}\le C\sqrt m
$
because the operator norm of $\mathcal M$ in $L^2(\mu)$ is bounded by some absolute constant.
The desired inequality follows by duality now. 

\section{The $L^2(\mu)$ bound for the Riesz transform of the $\Psi$-function}\label{l2bound}

This section is devoted to the proof of the inequality 
$$
\int_{\R^2} |R(\Psi m_2)|^2\,d\mu\le C_9 m\,.
$$
It will suffice to get a uniform bound of the same kind for each $\Psi\ci A$ ($A\ge 2)$ 
separately. We shall compare $R(\Psi\ci Am_2)$ to 
$\sum_j \chi\ci{\R^2\setminus 2AT_j}R(\chi\ci{\wt T_j}\mu)$. Note that the $L^2(\mu)$ 
norm of the latter is bounded by $C\sqrt m$, which can be shown by exactly the same duality
argument as in the previous section only with $(R^*)^\sharp$ instead of $\mathcal M$.

We start with estimating each difference 
$$
R\bigl(\tfrac{\mu(\wt T_j)}{\pi A^2r_j^2}\chi\ci{AT_j}m_2\bigr)-\chi\ci{\R^2\setminus 2AT_j}R(\chi\ci{\wt T_j}\mu)
$$
pointwise.
If $x\in 2AT_j$, then only the first term matters and we can use the trivial $L^\infty$ bound
$$
\bigl|R(\tfrac{\mu(\wt T_j)}{\pi A^2r_j^2}\chi\ci{AT_j}m_2)\bigr|\le C\frac{\mu(\wt T_j)}{(Ar_j)^s}\,.
$$
If $x\notin 2AT_j$, we can use the smoothness of the kernel $\frac{x-y}{|x-y|^{s+1}}$ and the 
cancellation property of the measure 
$\tfrac{\mu(\wt T_j)}{\pi A^2r_j^2}\chi\ci{AT_j}m_2-\chi\ci{\wt T_j}\mu$ to get the bound
$
C\frac{\mu(\wt T_j)}{|x-c_j|^s}\frac{Ar_j}{|x-c_j|}\,.
$

Combining these two bounds, we see that the difference under consideration is bounded by
$$
C\mu(\wt T_j)\sum_{A'\in \A, A'\ge A} \frac{A}{A'}\frac{1}{(A'r_j)^s}\chi\ci {A'T_j}\,,
$$
which implies that $R(\Psi\ci Am_2)$ differs from 
$\sum_j\chi\ci{\R^2\setminus 2AT_j}R(\chi\ci{\wt T_j}\mu)$
by at most $C\sum_{A'\in \A, A'\ge A} \frac{A}{A'}g\ci{A'}$. But the $L^2(\mu)$-norms
of the Marcinkiewicz functions $g\ci{A'}$ are uniformly bounded by $\sqrt{C_8 m}$.

\section{The bottom cover}\label{botcov}

Choose $\rho^*>0$ so small that the $\mu$-measure of the $\rho^*$-neighborhood of the union
of the boundaries of the top cover disks $T_j$ is less than $\e m$ and that $|R(\Psi m_2)|^2(x')-|R(\Psi m_2)|^2(x'')\le 1$ whenever $|x'-x''|\le 3\rho^*$ (note that $|R(\Psi m_2)|^2$
is a continuous function tending to $0$ at infinity).

Take any point $x\in \cup_j T_j$ whose distance to the boundary of any $T_j$ is greater than $\rho^*$
and choose some disk $D(x,t_0)$ with $0<t_0<\rho^*$ satisfying 
$$
\mu(D(x,Mt_0))\le \delta t_0^s\,.
$$
According to our assumptions, the points $x$ for which such disk does not exist form
a set of $\mu$-measure $0$. Now put $t_j=(1-3\e)^jt_0$ ($j\ge 1$). Let $k\ge 0$ be the least 
index such that
$$
\mu(D(x,t_k)\setminus D(x,t_{k+1}))\le 6\e \mu(D(x,t_k))\,.
$$
It may happen, of course, that this inequality never holds but in that case 
$$
\frac{\mu(D(x,t_{j+1}))}{t_{j+1}^2}\le \frac{1-6\e}{(1-3\e)^2}\frac{\mu(D(x,t_{j}))}{t_{j}^2}
$$
for all $j\ge 0$. Since $\frac{1-6\e}{(1-3\e)^2}<1$, this implies that at every such point
$x$, the measure $\mu$ has zero density with respect to $m_2$ whence such bad points form
a set of $\mu$-measure $0$.

Put $\rho(x)=t_k$. We claim that 
$$
\mu(D(x,M\rho(x)))\le (1-3\e)^{-s}M^s\delta \rho(x)^s\le 2M^s\delta \rho(x)^s,
$$
provided that $\e<0.01$, say.

If $M\rho(x)>t_0$, this follows from the choice of $t_0$ immediately. Otherwise, choose
the largest $j$ such that $t_j\ge M\rho(x)$. Note that the sequence 
$\frac{\mu(D(x,t_{i}))}{t_{i}^s}$ ($0\le i\le j$) is decreasing and its zeroth term is at most $\delta$.
Thus
$$
\mu(D(x,M\rho(x)))\le \delta t_j^s\le (1-3\e)^{-s}M^s\delta \rho(x)^s\,.
$$
Now use the Besicovitch covering lemma to find a finite sequence of disks $B_j=B(x_j,\rho(x_j))$
that has covering number not exceeding $C_9$ and covers all points outside an exceptional
set of measure at most $3\e m$ (which includes the points outside $\cup_j T_j$, the points
too close to the boundaries, various bad points, and a small extra piece that ensures that
the covering is finite rather than countable). We shall write $\rho_j$ instead of $\rho(x_j)$
from now on and assume that the sequence $\rho_j$ is non-increasing.

Let $\wt B_j=(1-3\e)B_j\setminus\cup_{i<j}B_i$. Note that the sets $\wt B_j$ cover all points
in the union $\cup_j B_j$ except those that lie in the set $\cup_j(B_j\setminus(1-3\e)B_j)$ 
whose $\mu$-measure does not exceed $6\e\sum_j\mu(B_j)\le 12C_9\e m$.

Another nice property of $\wt B_j$ is that the distance from $\wt B_i$ to $\wt B_j$ is at least
$3\e\max(\rho_i,\rho_j)$ (the ordering of $\rho_j$ was done exactly for this purpose). The sets
$\wt B_j$ are nice but they may be a bit too thin, so let us also introduce for each $j$ 
the set $\Omega_j$, which is the $\e\rho_j$-neighborhood of $\wt B_j$.
Ignoring the indices for which $\wt B_j=\varnothing$, we can say that the sets $\Omega_j$
are still well-separated: the distance from each $\Omega_j$ to any other $\Omega_i$ is at
least $\e\rho_j$, and each set $\Omega_j$ contains some disk $\wt \Omega_j$ of radius $\e\rho_j$.

The sets $\Omega_j$ will be used as the first generation Cantor cells. 

\section{The full $N$-level Cantor construction and the associated measure $\mu'$}\label{fulcan}

For the zeroth level, we put $Q_1^{(0)}=\R^2$, $\mu_1^{(0)}=\mu$, $H_1^{(0)}=H$, $m_1^{(0)}=m$.
 
For the first level, we put $Q_j^{(1)}=\Omega_j$, $\mu_j^{(1)}=\chi\ci{\wt B_j}\mu$, 
$H_j^{(1)}=2\mathcal H^s(\supp\mu_j^{(1)})$, $m_j^{(1)}=\mu_j^{(1)}(\R^2)/2=\mu(\wt B_j)/2$.

To get the second level, we repeat the entire construction for each measure $\mu_j^{(1)}$ instead of 
$\mu$ (using the corresponding parameters $m_j^{(1)}$ and $H_j^{(1)}$ instead of $m$ and $H$ but the 
same $\e,M,\delta$). We shall get some new cells $Q_j^{(2)}$. We can easily ensure that each 
$Q_j^{(2)}$ is contained in a unique cell $Q_i^{(1)}$ of the previous generation if we choose 
the radius bound $r^*$ (depending on $j$) for the top cover of $\mu_j^{(1)}$ small enough (note that $\supp \mu_j^{(1)}$ lies deep inside $Q_j^{(1)}$). It will be also convenient to assume that the radius bound $\rho^*$ for the bottom cover of $\mu_j^{(1)}$ is chosen so that $M\rho^*$ is much less than
all the distances from $Q_j^{(1)}$ to all other cells $Q_i^{(1)}$.

Continuing this procedure for $N$ steps, we get a Cantor structure $Q_j^{(n)}$ on the plane 
($n=0,\dots,N$). We define the rarefied measure $\mu'$ by
$$
\mu'=\sum_j\mu_j^{(N)}\,.
$$ 
Note that $\mu'$ is just the restriction of $\mu$ to some subset of the plane. 

The important points to keep in mind are the following:

\prop{Small measure loss} Since every time we go one level down we get only $C_{10}\e$-portion
of the entire measure outside the next level Cantor cells, we have
$$
\mu'(Q_j^{(n)})\ge(1-C_{10}\e)^{N-n}\cdot 2 m_j^{(n)}\ge m_j^{(n)}
$$ 
if we choose $\e$ so small that $(1-C_{10}\e)^{N}\ge\frac 12$. 

\prop{Subordination} $\mu'$ is dominated by $\mu_j^{(n)}$ on $Q_j^{(n)}$.

\prop{The total counts} For every fixed $n=0,\dots,N-1$, we have $\sum_j m_j^{(n)}\ge\frac m2$, 
$\sum_j H_j^{(n)}\le H$.

\section{Partial Riesz potentials $R^{(n)}\mu'$ and the key estimates}\label{parrie}

For every $x\in\supp\mu'$, denote by $Q^{(n)}(x)$ the unique set $Q_j^{(n)}$ containing $x$.
Put
$$
R^{(n)}\mu'(x)=\int_{Q^{(n)}(x)\setminus Q^{(n+1)}(x)}\frac{x-y}{|x-y|^{s+1}}\,d\mu'(y)\,,\qquad 
n=0,\dots,N-1\,.
$$

The key observation is that, once $N$ is fixed, the other three construction parameters 
$\e,M,\delta$ can be chosen so that the following three claims hold:

\prop{Claim 1} On $\supp\mu'$, one has
$$
\Bigl|\sum_{n=0}^{N-1}R^{(n)}\mu'\Bigr|\le R^\sharp\mu'+1\,.
$$ 
 
\prop{Claim 2} For every $n=0,\dots,N-2$, one has
$$
\Bigl|\int_{\R^2}\Bigl\langle
R^{(n)}\mu',\sum_{k=n+1}^{N-1}R^{(k)}\mu'
\Bigr\rangle\,d\mu'  \Bigr|\le \frac {m^{5/2}}{4N C_{20} H^2}\sum_{k=n+1}^{N-1}\|R^{(k)}\mu'\|\ci{L^2(\mu')}\,.
$$

\prop{Claim 3} 
$$
\int_{\R^2}|R^{(n)}\mu'|^2\,d\mu'\ge C_{20}^{-2}\frac{m^5}{H^4} 
$$
for all $n=0,\dots,N-1$.

Once these claims are established, we can finish the argument as follows.
On one hand, Claim 1 implies that
$$
\int_{\R^2}\Bigl|\sum_{n=0}^{N-1}R^{(n)}\mu'\Bigr|^2\,d\mu' 
\le \int_{\R^2}|R^\sharp\mu'+1|^2\,d\mu\le 2(C_4+1)^2 m\,.
$$
On the other hand, expanding the square and combining Claims 2 and 3, we get the lower
bound 
\begin{multline*}
\sum_{n=0}^{N-1}\|R^{(n)}\mu'\|\ci{L^2(\mu')}
\left(\|R^{(n)}\mu'\|\ci{L^2(\mu')}-\frac {m^{5/2}}{2 C_{20} H^2}\right)
\\
\ge
\frac 12\sum_{n=0}^{N-1}\|R^{(n)}\mu'\|^2\ci{L^2(\mu')}\ge\frac N{2C_{20}^2}\frac{m^5}{H^4}
\,.
\end{multline*}
If $N>4(C_4+1)^2C_{20}^2\left(\frac Hm\right)^4$, we get a clear contradiction.

\section{The proof of Claim 1}\label{claim1}

Let $x\in\supp\mu'$. 
Let $Q_j^{(N-1)}$ be the unique Cantor cell from the $(N-1)$-st level containing $Q^{(N)}(x)$.
Let $B$ be the disk in the bottom cover of $\mu_j^{(N-1)}$ that gave birth to $Q^{(N)}(x)$.
Let $\rho$ be its radius.
Recall that the radius bound $\rho^*$ in the construction of the bottom cover for the measure
$\mu_j^{(N-1)}$ was chosen much less than the distance from $Q_j^{(N-1)}$ to any other 
$Q_i^{(N-1)}$, so the disk $2B$ does not intersect any other $Q_i^{(N-1)}$. 
The value of the sum to estimate at the point $x$ can be written as 
$\int_{\R^2\setminus Q^{(N)}(x)}\frac{x-y}{|x-y|^{s+1}}\,d\mu'(y)$. It differs from the integral
over $\R^2\setminus 2B$ (which is dominated by $(R^\sharp\mu')(x)$ by the definition of the latter)
 only by the integral
$$
\int_{2B\setminus Q^{(N)}(x)}\frac{x-y}{|x-y|^{s+1}}\,d\mu'(y)\,.
$$
Now, the integrand is uniformly bounded by $\frac 1{(\e\rho)^s}$ and the measure is not greater
than $\mu_j^{(N-1)}(2B)\le \mu_j^{(N-1)}(MB)\le 2M^s\delta \rho^s$, provided that $M\ge 2$. Thus,
the integral is at most $1$, provided that $\frac{2M^s\delta}{\e^s}<1$.

\section{The oscillation bound}\label{oscbou}

Let $\nu$ be any finite (signed) measure. Assume that $\Omega\subset\R^2$ is 
contained in a disk $B=B(x,\rho)$ and is $\e\rho$-separated from the support of $\nu$.
Then
$$
\osc_\Omega R\nu\le \frac{2}{(\e\rho)^s}|\nu|(\tfrac M3B)+\frac{C_{11}}{M}
\sup_{r>0}\frac{|\nu|(D(x,r))}{r^s}\,.
$$ 
Indeed, take $x',x''\in\Omega$ and notice that the difference
$$
\frac{x'-y}{|x'-y|^{s+1}}-\frac{x''-y}{|x''-y|^{s+1}}
$$
is bounded by $\frac{2}{(\e\rho)^s}$ for all $y\in\supp\nu$ and by
$$
\frac{C\rho}{|x-y|^{s+1}}
$$
for $y\notin \frac M3B$ if $M\ge 6$, say. Integrating the first bound over $\frac M3B$ and the second
one over its complement with respect to $|\nu|$, we get the desired estimate.

We will also need the dual form of this estimate, which says that if $\nu$ is a finite positive
measure and $\eta$ is a signed measure supported on $\Omega$ with perfect cancellation ($\eta(\Omega)=0$), then
$$
\int_{\R^2} |R\eta|\,d\nu\le \left[\frac{2}{(\e\rho)^s}\nu(\tfrac M3B)+\frac{C_{11}}{M}
\sup_{r>0}\frac{\nu(D(x,r))}{r^s}\right]|\eta|(\Omega)\,.
$$ 

Similar bounds (with the same proofs, but, possibly, slightly
larger constants) hold for $R^*$ instead of $R$.

\section{Proof of Claim 2}\label{claim2}

Apply the obtained oscillation bound to $\Omega=Q_j^{(n+1)}\subset Q_i^{(n)}$ and the measure 
$\nu=\chi\ci{Q_i^{(n)}}\mu'$ which is dominated by $\mu_i^{(n)}$. Let $B$ be the disk 
in the bottom cover of $\mu_i^{(n)}$ that gave birth to the Cantor cell $Q_j^{(n+1)}$. Then
the first term in the oscillation bound does not exceed $\frac{4M^s\delta}{\e^s}$ and the second 
term is bounded by $\frac{C_3C_{11}}{M}$. Thus, for every Cantor cell $Q_j^{(n+1)}$ of the $(n+1)$-st generation, 
$$
\osc_{Q_j^{(n+1)}}R^{(n)}\mu'\le C_{12}\biggl(\frac{M^s\delta}{\e^s}+\frac 1M\biggr)\,.
$$  
On the other hand, the sum $\sum_{k=n+1}^{N-1}R^{(k)}\mu'$ has the cancellation property
\begin{multline*}
\int_{Q_j^{(n+1)}}\Bigl[\sum_{k=n+1}^{N-1}R^{(k)}\mu'\Bigr]\,d\mu'
\\
=
\iint_{x,y\in Q_j^{(n+1)}, Q^{(N)}(x)\ne  Q^{(N)}(y)}\frac{x-y}{|x-y|^{s+1}}\,d\mu'(x)\,d\mu'(y)=0\,.
\end{multline*}
Thus
\begin{multline*}
\Bigl|\int_{\R^2}\Bigl\langle
R^{(n)}\mu',\sum_{k=n+1}^{N-1}R^{(k)}\mu'
\Bigr\rangle\,d\mu'  \Bigr|\le 
C_{12}\left(\frac{M^s\delta}{\e^s}+\frac 1M\right)
\sum_{k=n+1}^{N-1}\|R^{(k)}\mu'\|\ci{L^1(\mu')}
\\
\le
C_{12}\left(\frac{M^s\delta}{\e^s}+\frac 1M\right)\sqrt{2m}\sum_{k=n+1}^{N-1}\|R^{(k)}\mu'\|\ci{L^2(\mu')}
\end{multline*}
by Cauchy-Schwarz, and Claim 2 will follow if $M$ and $\delta$
satisfy
$$
C_{12}\left(\frac{M^s\delta}{\e^s}+\frac 1M\right)\sqrt{2}\le \frac {m^{2}}{4N C_{20} H^2}\,.
$$

\section{The maximum principle}\label{maxpri}

Suppose that $\eta$ is a vector-valued measure with compactly supported $C^\infty$ 
density with respect to $m_2$. Then
$$
\max_{\R^2} R^*\eta=\max_{\supp\eta} R^*\eta
$$
provided that the left hand side is positive.

Indeed, the function $u=R^*\eta$ can be written as the $(s-1)$-dimensional 
Newton potential $U\nu$  where $\nu$ is some
scalar signed measure with compactly supported  $C^\infty$ density with respect to $m_2$ 
satisfying $\supp\nu\subset\supp\eta$ (see Section \ref{elepro}). We need the well-known fact that the density $p$ of $\nu$
can be recovered from the potential $u=U\nu$ by the formula
$$
p(x)=\sigma \int_{\R^2}\frac{u(x+y)-u(x)}{|y|^{5-s}}\,dm_2(y),
$$
where $\sigma$ is some non-zero real number and the integral is understood
in the principal value sense. 

To demonstrate it, define the class $S_\gamma$ ($\gamma>0$) of smooth functions in $\mathbb R^d$
by
$$
S_\gamma=\{\f\in C^\infty(\R^d):\f(x)=O(|x|^{-\gamma})\text{ as } |x|\to\infty\}.
$$
For $0<\text{Re}\,\alpha<\gamma$ and $\f\in S_\gamma$, define
$$
K_\alpha\f=A(d,\alpha)\,\f\ast\frac1{|x|^{d-\alpha}}\text{ \ with \ }
A(d,\alpha)=\pi^{\alpha-\frac{d}2}\frac{\Gamma(\frac{d-\alpha}2)}{\Gamma(\frac{\alpha}2)}\,.
$$
Note that for every $x\in\R^d$, $K_\alpha\f(x)$ is analytic in $\alpha$ in the strip $0<\text{Re}\,\alpha<\gamma$. The argument on pages 45--46 in \cite{Lan} shows that $K_\alpha\f(x)$ extends analytically to the wider strip  $-2<\text{Re}\,\alpha<\gamma$ and is given for  
$\text{Re}\,\alpha<0$ by the formula
$$
K_\alpha\f(x)=A(d,\alpha) \int_{\R^d}\frac{\f(x+y)-\f(x)}{|y|^{d-\alpha}}\,dy,
$$
where the integral converges absolutely for $\text{Re}\,\alpha>-1$ and should be understood as the principal value for  $-2<\text{Re}\,\alpha\le-1$. Note that Landkof writes $p$ instead of $d$ and $k_\alpha\ast$ instead of $K_\alpha$. Also, if $\f\in C_0^\infty(\R^d)$, then, for $\text{Re}\,\alpha\in(0,d)$, we have $K_\alpha\f\in S_{d-\text{Re}\,\alpha}$.

It is well-known that for $\f\in C_0^\infty(\R^d)$, we have
$$
K_\alpha K_\beta\,\f=K_{\alpha+\beta}\,\f,\quad\alpha,\beta>0,\quad\alpha+\beta<d.
$$
Also $K_0\f=\f$ (see \cite{Lan}, p.~46). Now, consider the identity
$K_\beta K_\alpha\,\f=K_{\beta+\alpha}\,\f$ for $0<\alpha<\min(2,d)$ and $0<\beta<d-\alpha$.
Viewing both parts of this identity as analytic functions of $\beta$ in the strip $0<\text{Re}\,\beta<d-\alpha$, we conclude that it holds in the entire strip and continues to hold in the wider strip
$-2<\text{Re}\,\beta<d-\alpha$ for the analytic extensions. Plugging $\beta=-\alpha$, we obtain
$K_{-\alpha} K_\alpha\,\f=\f$ for $0<\alpha<\min(2,d)$, which is equivalent to our reproduction formula for $d=2$, $\alpha=3-s$.

In particular, we can conclude that the integral on the right hand side vanishes for all 
$x\notin\supp\nu$. Now, since $u$ is smooth and tends to $0$ at infinity, the point of
maximum is guaranteed to exist if the maximum is positive. But then at the point of maximum, the
integral is certainly negative because the integrand is non-positive everywhere and negative 
for all sufficiently large $y\in\mathbb R^2$. Thus the point of maximum must belong to $\supp\nu\subset\supp\eta$, proving the claim.

We shall need a slightly more general fact below. If $\nu$ is a finite positive measure
with compactly supported  $C^\infty$ density with respect to $m_2$, and $g$ is any $C^\infty$
vector-valued function, then
$$
\max_{\R^2} [V(R\nu)+R^*(g\nu)]=\max_{\supp\nu} [V(R\nu)+R^*(g\nu)]\,,
$$
provided that the left hand side is positive.
 
Indeed, we can write $v(t)=\max_{\tau\ge 0}[\tau t-v^*(\tau)]$ where $v^*$ is the Legendre
transform of $v$ (all we really need to know is that $v^*\ge 0$). Thus,
$$
V(x)=\max_{\tau\ge 0,\,|e|=1} [\tau\langle e,x\rangle-v^*(\tau)]
$$  
and 
$$
V(R\nu)+R^*(g\nu)=\max_{\tau\ge 0,\,|e|=1}[ R^*((g-\tau e)\nu)-v^*(\tau)]\,.
$$ 
Again, if the maximum is positive, it is attained at some point $x$ and equals to
the value of  $R^*((g-\tau e)\nu)-v^*(\tau)$ at $x$ for some $\tau,e$. But then the maximum
of $R^*((g-\tau e)\nu)$ is also positive and is attained at some point $y\in\supp\nu$.
The chain of inequalities
\begin{multline*}
[V(R\nu)+R^*(g\nu)](y)\ge [R^*((g-\tau e)\nu)](y)-v^*(\tau)
\\
\ge 
[R^*((g-\tau e)\nu)](x)-v^*(\tau)=[V(R\nu)+R^*(g\nu)](x)
\end{multline*}
finishes the argument.

It will be convenient to restate the last result in the following form.
If $\Lambda>0$ and $V(R\nu)+R^*(g\nu)\le\Lambda$ on $\supp\nu$, then
$V(R\nu)+R^*(g\nu)\le\Lambda$ on the entire plane.

Note that this part fails dramatically for $s<1$ because the density
reproduction formula then becomes
more complicated and involves the Laplacian $\Delta u(x)$, which is (or, at least,
seems) totally out of control. 

\section{The mollified measure $\wt\mu$}\label{molmea}

We now return to the zeroth level of the Cantor structure and to the notation of Sections
\ref{topcov}--\ref{botcov}. For each disk $\wt\Omega_j$, choose some positive $C^\infty$
cap $\f_j$ such that $\supp\f_j\subset\wt\Omega_j$, 
$\|\f_j\|\ci{L^\infty(m_2)}\le \frac{\mu'(\Omega_j)}{(\e\rho_j)^2}$, and
$\int_{\R^2}\f_j\,dm_2=\mu'(\Omega_j)$. Put $\wt\mu_j=\f_j m_2$ and $\wt\mu=\sum_j\wt\mu_j$.

Our first task will be to get a decent growth bound for $\wt\mu$. Take any disk $D=D(x,r)$.
Write
$$
\wt\mu(D)=\sum_{j:\rho_j<r}\wt\mu(D\cap \Omega_j)+\wt\mu(D\cap(\cup_{j:\rho_j\ge r}\Omega_j))\,.
$$
Recall that $\Omega_j$ are disjoint (and even well-separated). Also note that every $\Omega_j$ 
with $\rho_j<r$ that intersects $D$ is contained in $3D$, whence the first sum does
not exceed 
$$
\sum_{j:\Omega_j\subset 3D}\wt\mu(\Omega_j)=
\sum_{j:\Omega_j\subset 3D}\mu'(\Omega_j)\le\mu'(3D)\,.
$$
On the other hand, on each $\Omega_j$ with $\rho_j\ge r$, the density of the measure $\wt\mu$ with
respect to $m_2$ is bounded by 
$$
\frac{\mu'(\Omega_j)}{(\e\rho_j)^2}\le \frac{\mu(B_j)}{(\e\rho_j)^2}
\le \frac{\mu(MB_j)}{(\e\rho_j)^2}\le 
\frac{2M^s\delta}{\e^2}\rho_j^{s-2}\le \frac{2M^s\delta}{\e^2}r^{s-2}\,, 
$$  
so the second term is at most $\frac{2\pi M^s\delta}{\e^2}r^s$. This yields the 
final growth bound
$$
\wt\mu(D)\le\mu'(3D)+\frac{2\pi M^s\delta}{\e^2}r^s\,,
$$
which can be used in two ways. First, choosing $\delta$ so that $\frac{2\pi M^s\delta}{\e^2}<1$,
we conclude that $\wt\mu(D)\le (3^sC_3+1)r^s=C_{13}r^s$ for all disks $D$. Second, taking $D=\frac M3B_j$, we 
conclude that
$$
\wt\mu(\tfrac M3B_j)\le 2M^s\delta\rho_j^s
+\frac{2\pi M^s\delta}{\e^2}\biggl(\frac{M}3\rho_j\biggr)^s
\le\frac{9M^{2s}\delta}{\e^2}\rho_j^s\,.
$$
We shall use these bounds in combination with the results of Section \ref{oscbou} in the next 
section. Now let us point out one more nice property of $\wt\mu$, which (in addition to having
an infinitely smooth density) is its great advantage over the unmollified measure $\mu'$: for 
every $j$,
$$
\|R(f\wt\mu_j)\|\ci{L^\infty(m_2)}\le 
C(\e\rho_j)^{2-s}\frac{\mu'(\Omega_j)}{(\e\rho_j)^2}\|f\|\ci{L^\infty(m_2)}\le
C_{14}\frac{M^s\delta}{\e^s}\|f\|\ci{L^\infty(m_2)}\,.
$$ 
The same bound holds for $R^*$ as well.

\section{The operator $\wt R$ and the mollified lower bound problem}\label{opewtR}

For a (signed) measure $\nu$ supported on $\cup_j\Omega_j$ and a point $x\in\Omega_j$, 
define $(\wt R\nu)(x)=(R(\chi\ci{\R^2\setminus\Omega(x)}\nu))(x)$ where $\Omega(x)$ is
the unique $\Omega_j$ containing $x$. Note that $\wt R\mu'=R^{(0)}\mu'$, of course.
The reason we introduce this new notation now is that we want to view $\wt R$ as an operator
while $R^{(0)}\mu'$ was rather a complex notation for a single function.

We want to compare $\int_{\R^2}V(\wt R\mu')\,d\mu'$ with $\int_{\R^2}V(R\wt \mu)\,d\wt\mu$ now. 
One remark about the notation may be in order. It would be slightly more accurate to say 
that the integrals are taken over $\cup_j\Omega_j$ because $\wt R\nu$ is defined only there.
Nevertheless, since we will integrate the expressions involving $\wt R$ exclusively with
respect to measures supported on $\cup_j\Omega_j$, we can view the integrals over $\R^2$ just as
integrals of functions defined almost everywhere rather than everywhere.

The comparison will be done in three steps.

\prop{Step 1} Since $V$ is Lipschitz with the Lipschitz constant $4$, we have
$$
\int_{\R^2}|V(\wt R\mu')-V(\wt R\wt\mu)|\,d\mu'\le 
4\int_{\R^2}|\wt R\mu'-\wt R\wt\mu|\,d\mu'\,.
$$

Let $\eta_j=\chi\ci{\Omega_j}\mu'-\wt\mu_j$. Note that $\wt R\eta_j=0$ on $\Omega_j$. 
Applying the dual form of the oscillation bound from Section \ref{oscbou} with 
$\eta=\eta_j$, $\nu=\chi\ci{\R^2\setminus\Omega_j}\mu'$, we get
$$
\int_{\R^2}|\wt R\eta_j|\,d\mu'\le 2\left(\frac{2M^s\delta}{\e^s}+\frac{C_{11}C_3}M\right)\mu'(\Omega_j)\,.
$$
Adding these estimates up, we conclude that 
$$
\int_{\R^2}|\wt R\mu'-\wt R\wt\mu|\,d\mu'\le C\left(\frac{M^s\delta}{\e^s}+\frac 1M\right)m
$$
and the same estimate holds for 
$\int_{\R^2}|V(\wt R\mu')-V(\wt R\wt\mu)|\,d\mu'$.

\prop{Step 2} The oscillation bound, combined with the growth bounds for $\wt\mu$ from
the previous section, implies that
$$
\osc_{\Omega_j} V(\wt R\wt\mu)\le 4\osc_{\Omega_j}\wt R\wt\mu
\le C\left(\frac{M^{2s}\delta}{\e^{2+s}}+\frac 1M\right)\,,
$$
so, since $\wt\mu(\Omega_j)=\mu'(\Omega_j)$ for all $j$, we have
$$
\Bigl|\int_{\R^2}V(\wt R\wt\mu)\,d\mu'-\int_{\R^2}V(\wt R\wt\mu)\,d\wt\mu\Bigr|\le
C\left(\frac{M^{2s}\delta}{\e^{2+s}}+\frac 1M\right)m\,.
$$
\prop{Step 3} Finally, recalling that $\|R\wt\mu_j\|\ci{L^\infty(m_2)}\le 
C_{14}\frac{M^s\delta}{\e^s}$ (see Section \ref{molmea}), we observe that
$$
\int_{\R^2}|V(\wt R\wt\mu)-V(R\wt\mu)|\,d\wt\mu\le 4C_{14}\frac{M^s\delta}{\e^s}m\,.
$$
Bringing all the above inequalities together, we obtain
$$
\int_{\R^2}V(\wt R\mu')\,d\mu'\ge \int_{\R^2}V(R\wt\mu)\,d\wt\mu-
C_{15}\left(\frac{M^{2s}\delta}{\e^{2+s}}+\frac 1M\right)m\,.
$$  
  
\section{The family of measures $\wt\mu^\alpha$ and the extremal problem}\label{fammea}

The direct estimate of $\int_{\R^2}V(R\wt\mu)\,d\wt\mu$ is still a hard task because,
despite we know that $V(R\wt\mu)$ has noticeable values on the plane, our maximum
principle, if we apply it to $V(R\wt\mu)$ directly, allows us only to conclude that 
$V(R\wt\mu)$ is not too small at some point on the support of $\wt\mu$, which seems 
next to useless for estimating any integral norm.

What saves the day is the idea of the equilibrium measure borrowed from the positive symmetric
kernel capacity theory.  Instead of proving the above energy type inequality for 
the original measure, we prove it for the ``energy minimizer'' $\wt\mu^a$ whose potential 
is, in some sense, almost constant on $\supp\wt\mu^a$, so an $L^\infty$ lower bound 
translates into an integral lower bound automatically.    The idea that the singular Riesz  potential of   extremal measure should be ``almost constant" (like in the classical potential theory with positive kernel) was somewhat explored in Section 5.2 of \cite{VoCBMS}. For $d=2, s=1$ the Cauchy potential was replaced by Menger's curvature potential which is again strange but positive kernel, see Tolsa's \cite{T2}.
  
To carry out the formal argument, consider all vectors $\alpha=\{\alpha_j\}$ with
non-negative entries and define $\wt\mu^\alpha=\sum_j\alpha_j\wt\mu_j$. Fix $\lambda>0$
and consider the functional
$$
\Phi(\alpha)=\lambda m\max_j \alpha_j+\int_{\R^2}V(R\wt\mu^\alpha)\,d\wt\mu^\alpha\,.
$$
Let $a$ be the minimizer of $\Phi(\alpha)$ under the constraint 
$\wt\mu^\alpha(\R^2)=\wt\mu(\R^2)$ (recall that $\wt\mu(\R^2)\in[m,2m]$). The minimizer exists because $\Phi(\alpha)$ is a continuous function of $\alpha$ tending
to $+\infty$ as $\max_j\alpha_j\to+\infty$.

Let us assume that $\int_{\R^2}V(R\wt\mu)\,d\wt\mu\le\lambda m$. Then $\Phi(a)\le 2\lambda m$ 
whence all $a_j\le 2$, so the extremal measure $\wt\mu^a$ is dominated by $2\wt\mu$.

Now let us fix any $j$ with $a_j>0$, take a small $t>0$, and try to replace $\wt\mu^a$ by $\left[1-t\wt\mu(\R^2)^{-1}\wt\mu_j(\R^2)\right]^{-1}(\wt\mu^a-t\wt\mu_j)$, which is also an 
admissible measure.

If we just subtract $t\wt\mu_j$ without the renormalization, $\max_ja_j$ will not increase and
the integral part will change in the first order by 
\begin{multline*}
-t\Bigl[\int_{\R^2}V(R\wt\mu^a)\,d\wt\mu_j+
\int_{\R^2}\langle \nabla V(R\wt\mu^a),R\wt\mu_j\rangle\,d\wt\mu^a\Bigr]
\\
=-t\int_{\R^2}[V(R\wt\mu^a)+
R^*(\nabla V(R\wt\mu^a)\wt\mu^a)]
\,d\wt\mu_j=-tI\,.
\end{multline*}
Since the renormalization can raise the value of any part of $\Phi(a)$ at most
$\left[1-t\wt\mu(\R^2)^{-1}\wt\mu_j(\R^2)\right]^{-3}$ times, we should have
$$
\left[1-t\wt\mu(\R^2)^{-1}\wt\mu_j(\R^2)\right]^{-3}(\Phi(a)-tI)\ge \Phi(a)-o(t)\qquad
\text{as }t\to 0+\,,
$$
whence 
$$
I\le 3\Phi(a)\wt\mu(\R^2)^{-1}\wt\mu_j(\R^2)\le 6\lambda\wt\mu_j(\R^2)
$$ 
because $\Phi(a)\le 2\lambda m$ and $\wt\mu(\R^2)\ge m$.

Thus, $V(R\wt\mu^a)+
R^*(\nabla V(R\wt\mu^a)\wt\mu^a)$ is at most $6\lambda$ on $\Omega_j$ on average (with
respect to the measure $\wt\mu_j$). Now notice that $|\nabla V|\le 4$ and $\wt\mu^a$ may have 
the growth bounds only twice worse than those for $\wt\mu$. The immediate conclusion is that 
$$
\osc_{\Omega_j}[V(R\wt\mu^a)+
R^*(\nabla V(R\wt\mu^a)\wt\mu^a)]\le C_{16}\left(\frac{M^{2s}\delta}{\e^{2+s}}+\frac 1M\right)
$$
(compare with Sections \ref{claim2} and \ref{opewtR}).

So 
$$
V(R\wt\mu^a)+
R^*(\nabla V(R\wt\mu^a)\wt\mu^a)\le 6\lambda+
C_{16}\left(\frac{M^{2s}\delta}{\e^{2+s}}+\frac 1M\right)=6\lambda+\beta
$$
on the entire $\Omega_j$ and, since $j$ was chosen arbitrarily, on $\supp\wt\mu^a$. But then this
estimate automatically extends to the entire plane by the maximum principle.

Integrating it against $\Psi dm_2$, we get
$$
(6\lambda+\beta)\int_{\R^2}\Psi\,dm_2\ge \int_{\R^2}V(R\wt\mu^a)\Psi\,dm_2+
\int_{\R^2}R^*(\nabla V(R\wt\mu^a)\wt\mu^a)\Psi\,dm_2\,.
$$

\section{Proof of Claim 3}\label{claim3}

Now it is time to bring up everything we know about the top cover and
the associated $\Psi$-function in one final effort. First, we have
seen in Section \ref{topcov} that 
$$
\int_{\R^2}\Psi\,dm_2\le C_6m\,.
$$ 
Second, the measure 
$\wt\mu^a$ satisfies the assumptions on the measure $\nu$ in Sections \ref{topcov}, \ref{functV}.
Thus
$$
\int_{\R^2}V(R\wt\mu^a)\Psi\,dm_2\ge C_7^{-1}\frac{m^3}{H^2}\,. 
$$
Third, the last remaining integral can be rewritten as 
$$
\int_{\R^2}\langle R(\Psi m_2),\nabla V(R\wt\mu^a)\rangle\,d\wt\mu^a
$$
which, by Cauchy-Schwarz, does not exceed
$$
\Bigl(\int_{\R^2}|R(\Psi m_2)|^2\,d\wt\mu^a\Bigr)^{1/2}
\Bigl(\int_{\R^2}|\nabla V(R\wt\mu^a)|^2\,d\wt\mu^a\Bigr)^{1/2}
$$
in absolute value.

Now, due to the second restriction on the radius bound $\rho^*$ in Section \ref{botcov}
and the inequality $\wt\mu^a\le 2\wt\mu$,
the first integral is bounded by
$$
2\int_{\R^2}|R(\Psi m_2)|^2\,d\mu'+2m\le 2(C_9+1)m
$$
according to the result of Section \ref{l2bound}. To estimate the second integral, we use
the inequality $|\nabla V|^2\le 4V$ from Section \ref{functV} and obtain
$$
\int_{\R^2}|\nabla V(R\wt\mu^a)|^2\,d\wt\mu^a
\le 4\int_{\R^2}V(R\wt\mu^a)\,d\wt\mu^a\le 4\Phi(a)\le 8\lambda m\,.
$$
Putting all these estimates together, we see that either $\lambda\le\beta$, or 
$$
7C_6\lambda\ge C_7^{-1}\left(\frac mH\right)^2-4\sqrt{C_9+1}\sqrt\lambda\,.
$$
Taking $\lambda=C_{17}^{-1}\left(\frac mH\right)^4$ with sufficiently
large $C_{17}$, and recalling that $m\le C_3H$
due to the growth bound, we see that the second possibility fails. So, either
our initial assumption $\int_{\R^2}V(R\wt\mu)\,d\wt\mu  \le \lambda m$ was false,
or $\lambda\le\beta$. In both cases, we conclude that
$$
\int_{\R^2}V(R\wt\mu)\,d\wt\mu  \ge \lambda m -\beta m=
C_{17}^{-1}\left(\frac mH\right)^4m-C_{16}\left(\frac{M^{2s}\delta}{\e^{2+s}}+\frac 1M\right)m\,.
$$
Recalling the comparison inequality between  $\int_{\R^2}V(\wt R\mu')\,d\mu'$ and  $\int_{\R^2}V(R\wt\mu)\,d\wt\mu$ from Section 19, we finally obtain
$$
\int_{\R^2}V(\wt R\mu')\,d\mu'  \ge
C_{17}^{-1}\left(\frac mH\right)^4m-C_{18}\left(\frac{M^{2s}\delta}{\e^{2+s}}+\frac 1M\right)m\,.
$$ 
Returning to the notation of Section \ref{fulcan} and considering $\mu_j^{(n)}$ instead of $\mu$,
we get the inequalities
$$
\int_{Q_j^{(n)}}V(R^{(n)}\mu')\,d\mu'  \ge
C_{17}^{-1}\left(\frac {m_j^{(n)}}{H_j^{(n)}}\right)^4m_j^{(n)}-C_{18}\left(\frac{M^{2s}\delta}{\e^{2+s}}+\frac 1M\right)
m_j^{(n)}\,.
$$ 
Summing these estimates over $j$ and taking into account that $V(x)\le|x|^2$, we 
arrive at the estimate
$$
\int_{\R^2}|R^{(n)}\mu'|^2\,d\mu'  \ge
C_{17}^{-1}\sum_j\left(\frac {m_j^{(n)}}{H_j^{(n)}}\right)^4m_j^{(n)}-C_{18}\left(\frac{M^{2s}\delta}{\e^{2+s}}+\frac 1M\right)
\sum_j m_j^{(n)}\,.
$$ 
According to Section 12 we have $2m\ge2\sum_jm_j^{(n)}\ge m$, which allows us to estimate the 
second sum from above by $m$. To estimate the first sum from below, note that for any
positive numbers $a_j, b_j$, we have
$$
\sum\frac{a_j^5}{b_j^4}\ge\frac{(\sum a_j)^5}{(\sum b_j)^4}\,,
$$
which is just the H\"older inequality
$$
\sum b_j^{4/5}\frac{a_j}{b_j^{4/5}}\le
\left[\sum(b_j^{4/5})^{5/4}\right]^{4/5}
\biggl[\sum\biggl(\frac{a_j}{b_j^{4/5}}\biggr)^5\biggr]^{1/5}
$$
in disguise. Applying it with $a_j=m_j^{(n)}$ and $b_j=H_j^{(n)}$,
we obtain the bound
$$
\int_{\R^2}| R^{(n)}\mu'|^2\,d\mu'  \ge
C_{19}^{-1}\left(\frac mH\right)^4m-C_{18}\left(\frac{M^{2s}\delta}{\e^{2+s}}+\frac 1M\right)m
$$ 
with $C_{19}=32C_{17}$.
Thus, we will get Claim 3 in Section 13 with $C_{20}=\sqrt{2C_{19}}$ if $M$ and $\delta$ 
satisfy 
$$
C_{18}\left(\frac{M^{2s}\delta}{\e^{2+s}}+\frac 1M\right)\le \frac 1{2C_{19}}\left(\frac mH\right)^4\,.
$$
It remains to note that, once $N$ and $\varepsilon$ are fixed,
 we can always choose first $M>1$ and then $\delta>0$ to 
satisfy this condition simultaneously with the conditions 
$$
\frac{2M^s\delta}{\varepsilon^s}<1 \quad \text{ and }\quad
C_{12}\left(\frac{M^s\delta}{\e^s}+\frac 1M\right)\sqrt{2}\le \frac {m^{2}}{4N C_{20} H^2}
$$ 
in Sections 14 and 16 correspondingly.

\section{Concluding remarks}

The same proof works in any dimension $d$ for $s\in(d-1,d)$. 
To cover the other values of $s$, we need some form of the maximum principle
(no matter how week; the equilibrium measure idea should allow one to turn any decent statement
of the kind ``small on the support, hence small everywhere'' into the desired
$L^2$ bound). Of course, more direct ways to get the lower bound may be even more interesting.

Notice also that we also proved the following theorem
\begin{theorem}
 Let $s\in(0,2)\setminus\{1\}$ and let $\mu$ be a strictly positive finite Borel measure in $\R^2$ such that $\mathcal H^s(\supp\mu)<+\infty$. Then $\sup_{\epsilon>0} |R_{\epsilon} \mu(x)|=\infty$ for  $\mu$ a. e. $x$,  and the operator norm
$\|R_\mu: L^2(\mu)\rightarrow L^2(\mu)\|=\infty$.
\end{theorem}
The same is again true for $s\in(d-1,d)$ in any $\R^d$. This is just reformulations of our main  theorem. This is easy to see by using several non-homogeneous Harmonic Analysis reductions as in \cite{NTV-Acta}, see also \cite{VoRMI}.

\end{document}